\documentclass[10pt,reqno]{amsart}
\usepackage{amsmath}
\usepackage{amssymb}
\usepackage{amsthm}

\newlength{\defbaselineskip}
\newcommand{\setlinespacing}[1]%
           {\setlength{\baselineskip}{#1 \defbaselineskip}}

\numberwithin{equation}{section}

\newtheorem{thm}{Theorem}[section]
\newtheorem{prop}[thm]{Proposition}
\newtheorem{lem}[thm]{Lemma}

\theoremstyle{definition}

\theoremstyle{remark}

\numberwithin{equation}{section}

\begin{document}

\title[Carleman inequalities and unique continuation]
{Carleman inequalities for fractional Laplacians and unique continuation}

\author{Ihyeok Seo}


\address{Department of Mathematics, Sungkyunkwan University, Suwon 440-746, Republic of Korea}
\email{ihseo@skku.edu}

\maketitle

\begin{abstract}
We obtain a unique continuation result for fractional Schr\"odinger operators with
potential in Morrey spaces.
This is based on Carleman inequalities for fractional Laplacians.
\end{abstract}

\section{Introduction}

The aim of this paper is to obtain a unique continuation result for the
fractional Schr\"odinger operator $(-\Delta)^{\alpha/2}+V(x)$, $0<\alpha<n$.
Recently, this operator has attracted interest from mathematics as well as mathematical physics.
This is because Laskin~\cite{La} introduced the fractional quantum mechanics governed by the fractional Schr\"odinger
equation
$$i\partial_t\Psi(x,t)=((-\Delta)^{\alpha/2}+V(x))\Psi(x,t),$$
where the fractional Schr\"odinger operator plays a central role.

More generally, we will consider the following differential inequality
\begin{equation}\label{equ}
|(-\Delta)^{\alpha/2}u(x)|\leq V(x)|u(x)|,\quad x\in\mathbb{R}^n,\quad n\geq3,
\end{equation}
where $(-\Delta)^{\alpha/2}$ is defined for $0<\alpha<n$
by means of the Fourier transform $\mathcal{F}f$ $(=\widehat{f}\,)$, as follows:
$$\mathcal{F}[(-\Delta)^{\alpha/2}f](\xi)=|\xi|^\alpha\widehat{f}(\xi).$$
The problem is now to find conditions on the potential $V(x)$ that imply
the unique continuation property which means that a solution of \eqref{equ}
vanishing in an open subset of $\mathbb{R}^n$ must vanish identically.

In the classical case $\alpha=2$, the property was extensively studied in connection
with the problem of absence of positive eigenvalues of the Schr\"odinger operator $-\Delta+V(x)$.
Among others, Jerison and Kenig \cite{JK} proved the property
for $V\in L_{\textrm{loc}}^{n/2}$, $n\geq3$.
Around the same time, an extension to $L_{\textrm{loc}}^{n/2,\infty}$ was obtained by Stein \cite{St3}
with the smallness assumption that
$$\sup_{a\in\mathbb{R}^n}\lim_{r\rightarrow0}\|\chi_{B(a,r)}V\|_{L^{n/2,\infty}}$$
is sufficiently small.
(Here, $\chi_{B(a,r)}$ denotes the characteristic function of the ball with center $a\in\mathbb{R}^n$
and radius $r>0$.)
Note that this assumption is trivially satisfied
for $V\in L_{\textrm{loc}}^{n/2}$ because $L_{\textrm{loc}}^{n/2}\subset L_{\textrm{loc}}^{n/2,\infty}$.
Also, the above-mentioned results later turn out to be optimal in the context of $L^p$ potentials (see \cite{KN,KT}).

Recently, there was an attempt \cite{S2} to deal with the fractional case where $n-1\leq\alpha<n$.
After that, the author \cite{S3} extended Stein's result completely to $0<\alpha<n$.
Namely, it turns out that \eqref{equ} has the unique continuation property for $V\in L_{\textrm{loc}}^{n/\alpha,\infty}$
with the corresponding smallness assumption that
$$\sup_{a\in\mathbb{R}^n}\lim_{r\rightarrow0}\|\chi_{B(a,r)}V\|_{L^{n/\alpha,\infty}}$$
is sufficiently small.
See also \cite{L,S} for higher orders where $\alpha/2$ are positive integers,
and for some fractional elliptic equations see \cite{FF,A,A2}.

In this paper we improve the class of potentials to
the Morrey class $\mathcal{L}^{\alpha,p}$ which is defined for $\alpha>0$ and $1\leq p\leq n/\alpha$ by
$$V\in\mathcal{L}^{\alpha,p}\quad\Leftrightarrow\quad
\|V\|_{\mathcal{L}^{\alpha,p}}
:=\sup_{Q\text{ cubes in }\mathbb{R}^{n}}|Q|^{\alpha/n} \Big( \frac{1}{|Q|} \int_{Q}V(y)^p dy \Big)^{\frac{1}{p}}<\infty.$$
In particular, $\mathcal{L}^{\alpha,p}=L^{p}$ when $p=n/\alpha$,
and even $L^{n/\alpha,\infty}\subset\mathcal{L}^{\alpha,p}$ for $p<n/\alpha$.
Our result is the following theorem.

\begin{thm}\label{thm}
Let $n\geq3$ and $0<\alpha<n$.
Assume that $V\in\mathcal{L}^{\alpha,p}$ for $p>(n-1)/\alpha$.
Let $u\in L^2\cap L^2(V)$ be a solution of \eqref{equ} vanishing
in a non-empty open subset of $\mathbb{R}^n$.
Then $u\equiv0$ if
\begin{equation}\label{small}
\sup_{a\in\mathbb{R}^n}\lim_{r\rightarrow0}\|\chi_{B(a,r)}V\|_{\mathcal{L}^{\alpha,p}}
\end{equation}
is sufficiently small. Here, $L^2(V)=L^2(V(x)dx)$.
\end{thm}

Let us give some remarks about the assumptions in the theorem.
First, $L^2\cap L^2(V)$ is the solution space for which we have unique continuation.
It should be noted that the space is dense in $L^2$.
In fact, consider $D_n=\{x\in\mathbb{R}^n:V^{1/2}\leq n\}$.
Then, for $f\in L^2$, $\chi_{D_n}f$ is contained in $L^2\cap L^2(V)$,
and $\chi_{D_n}f\rightarrow f$ as $n\rightarrow\infty$.
Now the Lebesgue dominated convergence theorem gives that $\chi_{D_n}f\rightarrow f$ in $L^2$.
Thus, the solution space is dense in $L^2$.

Next, by taking the rescaling
$u_\varepsilon(x)=u(\varepsilon x)$, the equation $(-\Delta)^{\alpha/2}u=Vu$
becomes $(-\Delta)^{\alpha/2}u_\varepsilon=V_\varepsilon u_\varepsilon$,
where $V_\varepsilon(x)=\varepsilon^\alpha V(\varepsilon x)$.
It is also easy to see that $\|V_\varepsilon\|_{\mathcal{L}^{\alpha,p}}=\|V\|_{\mathcal{L}^{\alpha,p}}$.
Hence, $\mathcal{L}^{\alpha,p}$ is invariant under the scaling.

\smallskip

The above theorem is a consequence of the following Carleman inequalities
which can be seen as natural extensions to the fractional Laplacians
$(-\Delta)^{\alpha/2}$ of those in \cite{CS} for the case $\alpha=2$.

\begin{prop}\label{prop}
Let $n\geq3$ and $0<\alpha<n$.
Assume that $V\in\mathcal{L}^{\alpha,p}$ for $p>(n-1)/\alpha$.
Then there exist sequence $\{t_m:m=0,1,...\}$ and constants $C,\beta>0$ independent of $m$ and $r$ such that
$$\big\|\chi_{B(0,r)}|x|^{-t_m-\frac{n-\alpha}{2}}f\big\|_{L^2(V)}
\leq C\|\chi_{B(0,r)}V\|_{\mathcal{L}^{\alpha,p}}^\beta
\big\||x|^{-t_m-\frac{n-\alpha}{2}}(-\Delta)^{\alpha/2}f\big\|_{L^2(V^{-1})}$$
for $f,(-\Delta)^{\alpha/2}f\in C_0^\infty(\mathbb{R}^n\setminus\{0\})$.
Here, $t_m\rightarrow\infty$ as $m\rightarrow\infty$.
\end{prop}

Throughout the paper, we will use the letter $C$ to denote positive constants possibly different at each occurrence.

\

\noindent\textbf{Acknowledgment.}
The author would like to thank Luis Escauriaza for bringing his attention to Carleman inequalities
for fractional Laplacians.

\section{Unique continuation}
Here we prove Theorem \ref{thm} assuming Proposition \ref{prop} which will be shown in the next section.

Without loss of generality, we may prove that the solution $u$ must vanish identically if
it vanishes in a sufficiently small neighborhood of the origin.

Since we are assuming that $u\in L^2\cap L^2(V)$ vanishes near the origin,
by \eqref{equ}, $(-\Delta)^{\alpha/2}u\in L^2(V^{-1})$ vanishes also near the origin.
Now, from the Carleman inequality in Proposition \ref{prop} (with a standard limiting argument involving a $C_0^\infty$ approximate identity),
one can easily see that
\begin{align}\label{lim}
\nonumber\big\|\chi_{B(0,r)}|x|^{-t_m-\frac{n-\alpha}{2}}&u\big\|_{L^2(V)}\\
&\leq C\|\chi_{B(0,r)}V\|_{\mathcal{L}^{\alpha,p}}^\beta
\big\||x|^{-t_m-\frac{n-\alpha}{2}}(-\Delta)^{\alpha/2}u\big\|_{L^2(V^{-1})}.
\end{align}
Note also that from \eqref{equ}
\begin{align*}
\big\||x|^{-t_m-\frac{n-\alpha}{2}}(-\Delta)^{\alpha/2}u\big\|_{L^2(V^{-1})}
&\leq C\big\|\chi_{B(0,r)}|x|^{-t_m-\frac{n-\alpha}{2}}u\big\|_{L^2(V)}\\
&+C\big\|(1-\chi_{B(0,r)})|x|^{-t_m-\frac{n-\alpha}{2}}(-\Delta)^{\alpha/2}u\big\|_{L^2(V^{-1})}.
\end{align*}
So, if we choose $r$ small enough so that $\|\chi_{B(0,r)}V\|_{\mathcal{L}^{\alpha,p}}^\beta$ is sufficiently small (see \eqref{small}),
then the first term on the right-hand side can be absorbed into the left-hand side of \eqref{lim}.
Thus we get
\begin{align*}
\big\|\chi_{B(0,r)}|x|^{-t_m-\frac{n-\alpha}{2}}u\big\|_{L^2(V)}
&\leq C\big\|(1-\chi_{B(0,r)})|x|^{-t_m-\frac{n-\alpha}{2}}(-\Delta)^{\alpha/2}u\big\|_{L^2(V^{-1})}\\
&\leq Cr^{-t_m-\frac{n-\alpha}{2}}\big\|(-\Delta)^{\alpha/2}u\big\|_{L^2(V^{-1})},
\end{align*}
which in turn implies
\begin{align*}
\big\|\chi_{B(0,r)}\big(\frac{r}{|x|}\big)^{t_m+\frac{n-\alpha}{2}}u\big\|_{L^2(V)}
\leq C\big\|(-\Delta)^{\alpha/2}u\big\|_{L^2(V^{-1})}<\infty.
\end{align*}
By letting $m\rightarrow\infty$ we conclude that $u=0$ on $B(0,r)$.
Now, $u\equiv0$ by a standard connectedness argument.

\section{Carleman inequalities}
In this section we will obtain the Carleman inequality in Proposition \ref{prop}
by using Stein's complex interpolation (\cite{St2}), as in \cite{CS},
on an analytic family of operators $S_z^{t,\alpha}$ defined by
$$S_z^{t,\alpha}g(x)=\frac{V(x)^{\frac z{2\alpha}}}{\Gamma((n-z)/2)}\int_{\mathbb{R}^n}K_z(x,y)V(y)^{\frac z{2\alpha}}g(y)dy,$$
where $0\leq \textrm{Re}\,z\leq n$ and
\begin{equation*}\label{ker}
K_z(x,y)=C_z\bigg(\frac{|y|}{|x|}\bigg)^{t+(n-z)/2}
\bigg(|x-y|^{-n+z}-\sum_{j=0}^{m-1}\frac1{j!}\Big(\frac{\partial}{\partial s}\Big)^j|sx-y|^{-n+z}\Big|_{s=0}\bigg).
\end{equation*}
Here, $S_z^{t,2}$ coincides with the analytic family of operators $S_z^t$ in \cite{CS}.
Note also that
\begin{equation}\label{22}
S_\alpha^{t,\alpha}\bigg(\frac{(-\Delta)^{\alpha/2}f(y)}{V(y)^{1/2}|y|^{t+(n-\alpha)/2}}\bigg)(x)
=\frac{f(x)V(x)^{1/2}}{|x|^{t+(n-\alpha)/2}}
\end{equation}
(see Lemma 2.1 in~\cite{S2}).

Let $m$ be nonnegative integers.
Now it is enough to show that there exist constants $C,\beta>0$ independent of $m$ and $r$ such that
\begin{equation}\label{St2}
\|\chi_{B(0,r)}S_{\alpha}^{t_m,\alpha} g\|_{L^2}\leq C\|\chi_{B(0,r)}V\|_{\mathcal{L}^{\alpha,p}}^{\beta}\|g\|_{L^2}
\end{equation}
for $q>(n-1)/\alpha$, $0<\varepsilon<\alpha(q-(n-1)/\alpha)$ and $t_m=m-1+(1-\varepsilon)/2$.
Indeed, from \eqref{22} and \eqref{St2},
\begin{equation*}
\bigg\|\chi_{B(0,r)}\frac{f(x)V(x)^{1/2}}{|x|^{t_m+(n-\alpha)/2}}\bigg\|_{L^2}
\leq C\|\chi_{B(0,r)}V\|_{\mathcal{L}^{\alpha,p}}^{\beta}
\bigg\|\frac{(-\Delta)^{\alpha/2}f(y)}{V(y)^{1/2}|y|^{t_m+(n-\alpha)/2}}\bigg\|_{L^2},
\end{equation*}
which is equivalent to
$$\bigg\|\chi_{B(0,r)}\frac{f(x)}{|x|^{t_m+(n-\alpha)/2}}\bigg\|_{L^2(V)}
\leq C\|\chi_{B(0,r)}V\|_{\mathcal{L}^{\alpha,p}}^{\beta}
\bigg\|\frac{(-\Delta)^{\alpha/2}f(y)}{|y|^{t_m+(n-\alpha)/2}}\bigg\|_{L^2(V^{-1})}$$
as desired.

To show \eqref{St2}, we use Stein's complex interpolation between the following two estimates
for the cases of $\textrm{Re}\,z=0$ and $n-1<\textrm{Re}\,z<\alpha q$:
\begin{equation}\label{JK}
\|\chi_{B(0,r)}S_{i\gamma}^{t_m,\alpha}g\|_{L^2}\leq Ce^{c|\gamma|}\|g\|_{L^2}
\end{equation}
and
\begin{equation}\label{St}
\|\chi_{B(0,r)}S_{n-1+\varepsilon+i\gamma}^{t_m,\alpha} g\|_{L^2}
\leq Ce^{c|\gamma|}\|\chi_{B(0,r)}V\|_{\mathcal{L}^{\alpha,p}}^{(n-1+\varepsilon)/2\alpha}\|g\|_{L^2},
\end{equation}
where $\gamma\in\mathbb{R}$, $p>(n-1)/\alpha$, $0<\varepsilon<\alpha(p-(n-1)/\alpha)$ and $t_m=m-1+(1-\varepsilon)/2$.
Indeed, since $n-1<n-1+\varepsilon<\alpha p\leq n$ and $p>1$,
we can easily get \eqref{St2} using the complex interpolation between \eqref{JK} and \eqref{St}.

It remains to show \eqref{JK} and \eqref{St}.
The first estimate \eqref{JK} follows immediately from Lemma 2.3 in \cite{JK}.
Indeed, consider the family of operators $T_z^t$ given by
$$T_z^tg(x)=\frac{1}{\Gamma((n-z)/2)}\int_{\mathbb{R}^n}H_z(x,y)g(y)|y|^{-n}dy,$$
where
$$H_z(x,y)=C_z|x|^{-t}|y|^{n+t-z}\bigg(|x-y|^{-n+z}-\sum_{j=0}^{m-1}\frac1{j!}\Big(\frac{\partial}{\partial s}\Big)^j|sx-y|^{-n+z}|_{s=0}\bigg).$$
Then it is clear that \eqref{JK} follows from
\begin{equation*}
\|T_{i\gamma}^{t_m}g\|_{L^2(dx/|x|^n)}\leq Ce^{c|\gamma|}\|g\|_{L^2(dx/|x|^n)}
\end{equation*}
which is Lemma 2.3 of \cite{JK}.

For the second one, we first recall from \cite{CS} (see (3.9) there) that
\begin{align*}
\bigg||x-y|^{-n+z}-\sum_{j=0}^{m-1}\frac1{j!}\Big(\frac{\partial}{\partial s}\Big)^j&|sx-y|^{-n+z}\Big|_{s=0}\bigg|\\
&\leq  Ce^{c|\textrm{Im}\,z|}\bigg(\frac{|x|}{|y|}\bigg)^{m-1+n-\textrm{Re}\,z}|x-y|^{-n+\textrm{Re}\,z}
\end{align*}
for $n-1<\textrm{Re}\,z<n$.
From this, we then get
$$|S_{n-1+\varepsilon+i\gamma}^{t_m,\alpha}g(x)|\leq Ce^{c|\gamma|}V(x)^{(n-1+\varepsilon)/2\alpha}
\int_{\mathbb{R}^n}|x-y|^{n-1+\varepsilon-n}V(y)^{(n-1+\varepsilon)/2\alpha}|g(y)|dy$$
if $0<\varepsilon<1$.
Hence it follows that
\begin{align}\label{end}
\nonumber\|\chi_{B(0,r)}&S_{n-1+\varepsilon+i\gamma}^{t_m,\alpha}g\|_{L^2}\\
&\leq Ce^{c|\gamma|}
\|\chi_{B(0,r)}I_{n-1+\varepsilon}(V(y)^{(n-1+\varepsilon)/2\alpha}|g(y)|)\|_{L^2(V^{(n-1+\varepsilon)/\alpha})},
\end{align}
where $I_\alpha$ denotes the fractional integral operator defined for $0<\alpha<n$ by
$$I_\alpha f(x)=\int_{\mathbb{R}^n}\frac{f(y)}{|x-y|^{n-\alpha}}dy.$$
Here we will use the following lemma to show \eqref{St},
which characterizes weighted $L^2$ inequalities for $I_\alpha$,
due to Kerman and Sawyer \cite{KeS} (see Theorem 2.3 there and also Lemma 2.1 in \cite{BBRV}):

\begin{lem}
Let $0<\alpha<n$. Assume that $w$ is a nonnegative measurable function on $\mathbb{R}^n$.
Then there exists a constant $C_w$ depending on $w$ such that
the following two equivalent estimates
\begin{equation*}
\|I_{\alpha/2}f\|_{L^2(w)}\leq C_w\|f\|_{L^2}
\end{equation*}
and
$$\|I_{\alpha/2}f\|_{L^2}\leq C_w\|f\|_{L^2(w^{-1})}$$
are valid for all measurable functions $f$
if and only if
\begin{equation}\label{ksc}
\sup_{Q}\bigg(\int_Qw(x)dx\bigg)^{-1}\int_Q\int_Q\frac{w(x)w(y)}{|x-y|^{n-\alpha}}dxdy
\end{equation}
is finite.
Here the sup is taken over all dyadic cubes $Q$ in $\mathbb{R}^n$,
and the constant $C_w$ may be taken to be a constant multiple of the square root of \eqref{ksc}.
\end{lem}

Indeed, it is known that $\|w\|_{\mathcal{L}^{\alpha,p}}<\infty$ for $p>1$ is a sufficient condition for
the finiteness of \eqref{ksc} (see Subsection 2.2 in \cite{BBRV}).
Namely, $\eqref{ksc}\leq C\|w\|_{\mathcal{L}^{\alpha,p}}$ for $p>1$.
Using this fact and applying the above lemma with $\alpha=n-1+\varepsilon$,
from \eqref{end} and $I_{\alpha/2}I_{\alpha/2}=I_{\alpha}$, we see that
for $1<q\leq n/(n-1+\varepsilon)$
\begin{align*}
\|\chi_{B(0,r)}&S_{n-1+\varepsilon+i\gamma}^{t_m,\alpha}g\|_{L^2}\\
&\leq Ce^{c|\gamma|}\|\chi_{B(0,r)}V^{(n-1+\varepsilon)/\alpha}\|_{\mathcal{L}^{n-1+\varepsilon,q}}^{1/2}
\|I_{(n-1+\varepsilon)/2}(V(y)^{(n-1+\varepsilon)/2\alpha}|g(y)|)\|_{L^2}\\
&\leq Ce^{c|\gamma|}\|\chi_{B(0,r)}V^{(n-1+\varepsilon)/\alpha}\|_{\mathcal{L}^{n-1+\varepsilon,q}}^{1/2}\\
&\qquad\qquad\qquad\qquad\times\|V^{(n-1+\varepsilon)/\alpha}\|_{\mathcal{L}^{n-1+\varepsilon,q}}^{1/2}
\|V^{(n-1+\varepsilon)/2\alpha}g\|_{L^2(V^{-(n-1+\varepsilon)/\alpha})}\\
&=Ce^{c|\gamma|}\|\chi_{B(0,r)}V^{(n-1+\varepsilon)/\alpha}\|_{\mathcal{L}^{n-1+\varepsilon,q}}^{1/2}
\|V^{(n-1+\varepsilon)/\alpha}\|_{\mathcal{L}^{n-1+\varepsilon,q}}^{1/2}\|g\|_{L^2}\\
&=Ce^{c|\gamma|}\|\chi_{B(0,r)}V\|_{\mathcal{L}^{\alpha,q(n-1+\varepsilon)/\alpha}}^{(n-1+\varepsilon)/2\alpha}
\|V\|_{\mathcal{L}^{\alpha,q(n-1+\varepsilon)/\alpha}}^{(n-1+\varepsilon)/2\alpha}\|g\|_{L^2}.
\end{align*}
Since $(n-1)/\alpha<(n-1+\varepsilon)/\alpha<q(n-1+\varepsilon)/\alpha\leq n/\alpha$
and $V\in\mathcal{L}^{\alpha,p}$,
by choosing $q,\varepsilon$ so that $p=q(n-1+\varepsilon)/\alpha$,
we now get the desired estimate \eqref{St}.


\end{document}